\documentclass[12pt]{amsart}

\usepackage{epsfig,color}

\usepackage{graphicx,psfrag}

\newtheorem{prop}{Proposition}
\newtheorem{defin}{Definition}
\newtheorem{thm}{Theorem}

\newtheorem{corol}{Corollary}
\newtheorem{lem}{Lemma}

\newtheorem{rem}{Remark}

\newcommand{\cc}{{\mathfrak C}}

\newcommand{\HH}{{\mathbb H}^2}

\newcommand{\N}{{\mathbb N}}

\newcommand{\PP}{{\mathcal P}}

\newcommand{\R}{{\mathbb R}}\newcommand{\RR}{{\mathbb R}^2}
\newcommand{\unisph}{{\mathbb S}^2}

\newcommand{\Z}{{\mathbb Z}}

\newcommand{\ar}{\text{area}}  %%%%%%%%%% AREA

\newcommand{\bo}{\partial} %%%%%%%%%% BOUNDARY
\newcommand{\ex}{\text{ex}}
\newcommand{\hyp}{\HH} %%%%%%%%%%%%% HYPERBOLIC PLANE
\newcommand{\iso}{\text{Iso}}   %%%%%%%%%%%%%%%% ISOMETRIES
\newcommand{\op}{\text{opt}}   %%%%%%%%%%%%% OPTICAL
\newcommand{\sph}{S^2} %%%%%%%%%%%%% SPHERE
\newcommand{\val}{\text{val}}  %%%%%%%%%% VALENCE
\newcommand{\vol}{\text{vol}}  %%%%%%%%%% VOLUME

\newcommand{\al}{\alpha}
\newcommand{\be}{\beta}
\newcommand{\ga}{\gamma}\newcommand{\Ga}{\Gamma}

\newcommand{\ep}{\varepsilon}

\newcommand{\ka}{\kappa}
\newcommand{\sig}{\sigma}
\newcommand{\tht}{\theta}\newcommand{\Tht}{\Theta}
\newcommand{\om}{\omega}
\newcommand{\vp}{\varphi}

\begin{document}

\bibliographystyle{plain}

\title[Complexity, etc]
{Growth rates  for geometric complexities and counting functions
in polygonal billiards}

\author{Eugene Gutkin and Michal Rams}
\address{IMPA, Rio de Janeiro, Brasil and UMK, Torun,
Poland;\hfill\hfill IMPAN, Warszawa, Poland}
\email{gutkin@impa.br,gutkin@mat.uni.torun.pl;rams@impan.gov.pl}

\keywords{Geodesic polygon, billiard map, billiard flow,
complexity, counting functions, unfolding of orbits, covering
space, exponential map}

%\subjclass{30F35, 30F60, 37B05, 37E35}

%\date{June 13, 2006}
\date{\today}

\begin{abstract}
We introduce a new method for estimating the growth of various
quantities arising  in dynamical systems. We apply our method to
polygonal billiards on surfaces of constant curvature. For
instance, we obtain power bounds of degree two plus epsilon for
%the number of
billiard orbits between almost all pairs of points in a
planar polygon.
\end{abstract}

\maketitle

\section*{Introduction and overview}   \label{intro}
Complexity of a dynamical system is measured with respect to a
coding of its orbits. The coding, in turn, is determined by
partitioning the phase space of the system into elementary pieces.
For dynamical systems with singularities, such as polygonal
billiards, connected components in the complement to the singular
set yield a natural partition. Convexity of its atoms
%(i. e., maximal domains of regularity for the billiard map)
with respect
to the geodesic structure in the phase space imposed by geometric
optics, is crucial in the study of billiard complexity
\cite{GT04}.

In the present study, as well as in \cite{GT04}, $P$ is a geodesic
polygon in a surface of constant curvature. Let, for concreteness,
$P$ be a planar polygon. We denote by $f_P(n)$ the number of words
of length $n$ generated by coding billiard orbits by visited
domains of regularity. When $P$ is simply connected, this
coincides with the coding by sides in $P$. It is known that
$f_P(n)$ is subexponential in $n$ \cite{GKT,GH97}, and for general
$P$ no better bound is known. If $P$ is a {\em rational polygon}
(i. e., its angles are commensurable with $\pi$ \cite{Gut86}),
$f_P(n)=O(n^3)$ \cite{CHT,GT04}. The current conjecture is that
for any planar polygon $f_P(n)=O(n^d)$ \cite{Gut03}.

In order to advance the understanding of billiard complexity, we
introduce the notion of {\em partial complexities}. Let $\Psi$ be
the phase space, and let $\PP$ be the defining partition.
Iterating the dynamics we obtain an increasing tower $\PP_n$ of
partitions; the full complexity is $f(n)=|\PP_n|$. If
$R\subset\Phi$, let $\PP_n(R)$ be the induced tower of its
partitions. The partial complexity based on $R$ is
$f_R(n)=|\PP_n(R)|$. Particular partial complexities have been
studied earlier. For instance, in \cite{GT96} we obtained
polynomial bounds on {\em direction complexity}, which is one of
the partial complexities investigated here.

In this work we introduce a new general approach to estimating
partial complexities. The setting is as follows. There is a family
of subsets $R_{\theta}$ foliating the phase space. Let
$f_{\theta}(n)$ be the partial complexity with base $R_{\theta}$.
Let $g_{\theta}(n)$ be the {\em counting function for singular
billiard orbits} starting from $R_{\theta}$. Under appropriate
assumptions, $f_{\theta}(n)$ and $g_{\theta}(n)$ have the same
growth, as $n\to\infty$. See section~\ref{comp_count}.

Let $\theta\in\Theta$, the parameter space. Suppose that we bound
the {\em average counting function}
$G(n)=\int_{\Theta}g_{\theta}(n)$. Tchebysheff inequality and the
zero-one law yield bounds for individual $g_{\theta}(n)$ valid for
almost all $\theta\in\Theta$. See section~\ref{bound}. Combined
with preceding remarks, these yield estimates on partial
complexities for almost all values of the parameter.

This is the general scheme for our approach to partial
complexities. This work implements this scheme for polygonal
billiards. We will now describe the contents of the paper in more
detail.

\medskip

In section~\ref{count} we investigate counting functions and their
averages. We establish the relevant framework in sufficient
generality, with the view towards a broad range of
geometric-dynamic  applications. The main results are
Propositions~\ref{crof_disc_prop} and~\ref{crof_cont_prop}
respectively. These yield geometric formulas for averages of
counting functions which are valid under mild assumptions of
transversality type.

Section~\ref{bound} is analytic, and also quite general. The
setting is as follows. There is a family of positive functions,
$g_{\theta}(p)$, of positive argument ($p\in\N$ and $p\in\R_+$ in
the discrete and continuous cases respectively), depending on
parameter $\theta\in\Theta$. Set
$G(p)=\int_{\Theta}g_{\theta}(p)d\theta$. From upper bounds on
$G(p)$ we derive estimates on individual $g_{\theta}(p)$; they are
valid for almost all $\theta$. Precise formulations depend on the
details of the situation. See Propositions~\ref{appli_prop1}
and~\ref{appli_prop2}.

Section~\ref{exa} sets the stage for applications to billiard
dynamics. Our billiard table is a geodesic polygon, $P$, in a
simply connected surface of constant curvature.
%We refer the
%reader to \cite{KH} for the basics on billiard dynamics.
%We relegate more specialized information on polygonal billiards to
%the appendix section~\ref{cover}.
There are two versions of billiard dynamics: the {\em billiard
flow} and the {\em billiard map}. In our discussion of partial
complexities, it is convenient to treat them separately.
Accordingly, section~\ref{exa} consists of several subsections;
each subsection deals with a particular partial complexity for a
particular kind of billiard dynamics.

We use two geometric parameters for partial complexities: the
direction and the position. The {\em direction complexity} tells
us how the set of phase points starting in the same direction
splits after bouncing off of the sides of $P$. The direction
complexity is defined for planar polygons. The {\em position
complexity} tells us about the splitting of beams of billiard
orbits emanating from a point of $P$. It is defined in all cases.

In each of the subsections of section~\ref{exa} we define a
counting function and check the assumptions of
section~\ref{count}; then we evaluate the integral over the
parameter space, i. e., we compute the average counting functions.
It turns out that they have geometric meanings. Here is a sample
of results from section~\ref{exa}. Let $G_P(l)$ be the average
position counting function for the billiard flow in a geodesic
polygon $P$. For planar polygons we have $G_P(l)=c_0(P)l^2$. See
Corollary~\ref{pos_count_flow_cor} in
section~\ref{pos_count_flow}. For polygons in $\unisph$ we have
$G_P(l)=c_+(P)l+c_+'(P)f(l)$ where $f$ is a universal periodic
function. See Corollary~\ref{sph_pos_count_cor} in
section~\ref{sph_pos_count_flow}. For polygons in $\hyp$ we have
$G_P(l)=c_-(P)\cosh l$. See Corollary~\ref{hyp_pos_count_cor} in
section~\ref{hyp_pos_count_flow}. The coefficients in these
formulas depend on how many corners $P$ has and on the number of
obstacles in its interior.

%\medskip

Section~\ref{comp_count}, again, is quite general. In this section
we obtain relationships between partial complexities with
one-dimensional base sets and counting functions. The main result
of this section is Proposition~\ref{comp_count_prop}. It says that
if the bases are one-dimensional, then the difference between the
partial complexity and the counting function is bounded, as time
goes to infinity. Other assumptions on the base have to do with
convexity in the phase space. The framework of this section is
that of {\em piecewise convex transformations} \cite{GT04}.

In section~\ref{appli} we specialize again to polygonal billiards.
Combining the material of preceding sections, we obtain bounds on
the position and direction complexities for the billard flow and
the billiard map. Here is a sample of our results. Let $P$ be a
euclidean polygon. Let $\tht\in S^1$ (resp. $z\in P$) be any
direction (resp. position). Let $fd_{\theta}(n)$ (resp. $h_z(l)$)
be the direction complexity for the billiard map (resp. position
complexity for the billiard flow). Then for almost all directions
$\theta$ (resp. for almost all positions $z$) we have
$fd_{\theta}(n)=O(n^{1+\ep})$ (resp. $h_z(l)=O(l^{2+\ep})$), where
$\ep>0$ is arbitrary. See Corollary~\ref{dir_est_cor2} and
Corollary~\ref{quad_log_cor}. Let now $P$ be a spherical polygon,
and let $h_z(l)$ be the position complexity for the billiard flow
in $P$. Then for almost every $z\in P$ there is a $C=C(z)$ and
arbitrarily large $l$ such that $h_z(l)\le Cl$. See
Corollary~\ref{sph_est_cor1}. For any $\ep>0$ and almost every
$z\in P$ we have $h_z(l)=O(l^{1+\ep})$. See
Corollary~\ref{sph_est_cor2}.

\vspace{2mm}

In the study of polygonal billiards the device of {\em unfolding
billiard orbits} is indispensable \cite{Gut86}. If $P\subset M$,
and $\be$ is a billiard orbit in $P$, its unfolding is a geodesic
in $M$. Several arguments in section~\ref{exa} use the technique
of lifting billiard orbits to the {\em universal covering space}
of $P$.\footnote{Not to be confused with the concept of universal
covering space in topology.} This notion was not written up in the
billiard literature. In our Appendix section~\ref{cover} we
present the relevant definitions and propositions.
Proposition~\ref{straigh_prop} puts forward the main property of
the universal covering space of a geodesic polygon. It relates the
unfoldings and the liftings of billiard orbits. The proofs in
section~\ref{exa} use Corollary~\ref{pull_back_cor} of
Proposition~\ref{straigh_prop}, which deals with the pullbacks of
lebesgue measures under unfoldings.

\medskip

In order to put  our results into perspective, we will now briefly
survey the literature on billiard complexities. The subexponential
growth of (full) billiard complexity for arbitrary euclidean
polygons is established in \cite{GKT} and \cite{GH97}. Both proofs
are indirect, in that they do not yield explicit subexponential
bounds. On the other hand, for rational  euclidean polygons the
complexity is cubic. This is contained in \cite{CHT} for convex
and in \cite{GT04} for all rational polygons. The arguments in
\cite{CHT} and \cite{GT04} rely on a theorem in \cite{M90}; it
says that the number of billiard orbits between any pair of
corners in a rational polygon grows quadratically in length. From
our viewpoint, this is a statement about the position counting
functions $g_z(l)$. It says that $g_z(l)=O(l^2)$ if $P\subset\RR$
is rational and $z\in P$ is a corner. By comparison, our
Corollary~\ref{quad_log_cor} and Proposition~\ref{comp_count_prop}
yield that $g_z(l)=O(l^{2+\ep})$ for any $\ep>0$ and almost all
$z\in P$ where $P\subset\RR$ is an arbitrary polygon. The
directional complexity $fd_{\tht}(n)$ has been studied in
\cite{GT96} and \cite{Hub95}. The work \cite{Hub95} concerns the
directional complexity for the billiard in a rational, planar
polygon $P$. Assume that $P$ is convex. Then \cite{Hub95} derives
an explicit formula for $fd_{\tht}(n)$, valid for {\em minimal
directions} $\tht$. (The set of nonminimal directions is
countable.) By this formula, $fd_{\tht}(n)=O(n)$. On the other
hand, \cite{GT96} shows that $fd_{\tht}(n)=O(n^d)$ for any
$P\subset\RR$ and an arbitrary $\tht$. The degree $d$ in the bound
does not depend on $\tht$. Our Corollary~\ref{dir_est_cor2}
estimates the complexity $fd_{\tht}(n)$ for an arbitrary polygon
$P\subset\RR$. It says that $fd_{\tht}(n)=O(n^{1+\ep})$ for any
$\ep>0$ and almost all directions $\tht$.

\vspace{2mm}

It is plausible that the bounds like Corollary~\ref{quad_log_cor},
Corollary~\ref{dir_est_cor2}, etc hold for any point $z\in P$, any
direction  $\tht\in S^1$, etc.

%\newpage

\section{Averages of counting functions}   \label{count}
In this section we introduce the framework of  counting functions
in differentiable dynamics. We will apply it to the billiard
dynamics later on. Our phase spaces are ``manifolds". By this we
will mean compact manifolds with boundaries, corners, and singular
points, in general. Our  setting involves i) a foliation of the
phase space by closed submanifolds that are fibers for a
projection onto a manifold of smaller dimension; ii) a submanifold
in the phase space, transversal to the fibers; iii) a {\em weight
function} on the product of the phase space and the time. See the
details below. The dynamics in question may be discrete or
continuous. We will expose the two cases separately. The two
subsections that follow are parallel, and the treatments differ in
technical details.

\subsection{Discrete dynamics}  \label{count_disc_sub}
Let  $T:X\to X,T^{-1}:X\to X$ be piecewise diffeomorphisms with
the following data.

\noindent 1. There is a fibration $\eta:X\to\Tht$ whose base is a
compact manifold and whose fibers $R_{\tht}=\eta^{-1}(\tht)$ are
compact submanifolds, such that
$\dim(R_{\tht})=\dim(X)-\dim(\Tht)$. We will use the notation
$X=\cup_{\tht\in\Tht}R_{\tht}$.

\noindent 2. There is a closed submanifold, $Y\subset X$,
$\dim(Y)=\dim(\Tht)$, such that for $k\in-\N$\footnote{By
convention, $\N=0,1,\dots$.} the manifolds $T^k(Y)$ are
transversal to the fibers $R_{\tht}$.

\noindent 3. There is a {\em weight function}, i. e., a
continuous, non-negative function $w(x,t)$ on $X\times\N$. The
function $w$ may depend only on time, e. g., $w=\chi_n$, the
indicator function of $[0,n-1]$.

\begin{rem}  \label{except_rem}
{\em Condition 2 may be weakened, as follows.

\vspace{1mm}

\noindent 2$'$. There is a closed submanifold, $Y\subset X$, and a
set $\Tht_{\ex}\subset\Tht$ of measure zero such that for
$k\in-\N$ and $\tht\in\Tht\setminus\Tht_{\ex}$ the manifolds
$T^k(Y)$ and $R_{\tht}$ are transversal. All of our results remain
valid if we replace condition 2 by the weaker condition 2$'$.
However, in our applications to polygonal billards, condition 2
may not hold only for polygons in surfaces of positive curvature.
See section~\ref{sph_pos_count_flow}. To simplify the exposition,
we will assume in what follows that $\Tht_{\ex}=\emptyset$.
%We leave details to the reader.
}
\end{rem}

%
%\medskip

In view of condition 2, $\Ga(\tht)=\{(x,k):\,x\in
R_{\tht},k\in\N,T^k(x)\in Y\}$ is a countable (at most) set. The
sets $\Ga_k(\tht)=\{(x,k):\,x\in R_{\tht},T^k(x)\in Y\}$ are
finite for all $k\in\N$, and $\Ga(\tht)=\cup\Ga_k(\tht)$.

%\medskip

We define the {\em weighted counting function} by
\begin{equation}            \label{coun_disc_eq}
g(\tht;w)=\sum_{(x,k)\in\Ga(\tht)}w(x,k).
\end{equation}
The {\em pure counting function} $g_n(\tht)$ corresponds to the weight $w=\chi_n$.
We have
\begin{equation}      \label{pur_coun_disc_eq}
g_n(\tht)= \sum_{k=0}^n|\Ga_k(\tht)|.
\end{equation}

\medskip

\begin{prop}   \label{crof_disc_prop}
Let $d\tht,dy$ be finite, lebesgue-class measures on $\Tht,Y$
respectively. Then for $k\in\N$ there are functions $r_k(\cdot)\ge
0$ on $Y$, determined by the data 1) and 2) alone, such that
\begin{equation}            \label{crof_disc_eq}
\int_{\Tht}g(\tht;w)d\tht\ =\ \int_{Y}\left\{\sum_{k\in\N}w(T^{-k}\cdot y,k)r_k(y)\right\}dy.
\end{equation}
%%
%The functions $r_k(\cdot)$  and the measure $dy$ are determined by the data 1) and 2)
%above; they do not depend on the weight $w$.
\begin{proof}
For any $k\in\N$ set $f_k=\eta\circ T^{-k}:Y\to\Tht$. By
conditions 1 and 2, $f_k$ is a local diffeomorphism. Therefore
$f_k^*(d\tht)=r_k(y)\,dy$.

It suffices to establish equation~\eqref{crof_disc_eq} for the
special case $w(x,i)=0$ if $i\ne k$. A point $x\in X$ contributes
to the integral in the left hand side of
equation~\eqref{crof_disc_eq} iff $T^k\cdot x\in Y$, or
equivalently, $\eta(x)=f_k(y),\,y\in Y$. The claim follows by a
straightforward change of variables.
\end{proof}
\end{prop}

%
%\begin{defin}    \label{reg_def}
%The transversal submanifold $Y$ is {\em regular} if all measured $d_ky$ coincide.
%\end{defin}
%

%\newpage

\subsection{Continuous dynamics}  \label{count_cont_sub}
Let $b^t:\Psi\to\Psi$ be a flow of piecewise diffeomorphisms on a
phase space $\Psi$ with the following data.

\noindent 1. There is a fibration $q:\Psi\to Z$ with a compact
base  and fibers $q^{-1}(z)=R_z\subset\Psi$, transversal to the
flow. We will use the notation $\Psi=\cup_{z\in Z}R_z$.

%\medskip

\noindent 2. There is a closed submanifold, $M\subset\Psi$,
$\dim(M)=\dim(Z)-1$, transversal to the flow, and such that
$N=\cup_{t\in\R_-}b^t\cdot M$ is transversal to the fibers
$R_z$.\footnote{Our results remain valid if the set of parameters
$Z_{\ex}\subset Z$ where the transversality fails has measure
zero. See Remark~\ref{except_rem}. In what follows, by condition
2$'$ we will mean the weakened condition 2 either in the setting
of section~\ref{count_cont_sub} or section~\ref{count_disc_sub}.}

\noindent 3. There is a  {\em weight function}, i. e., a
continuous, non-negative function $w(x,t)$ on $\Psi\times\R_+$. In
a special case, $w$ depends only on time, e. g., $w=\chi_l$, the
indicator function of $[0,l]$.

\medskip

In view of condition 2, $G(z)=\{(x,t):\,x\in R_z,0 \le t,b^t(x)\in
M\}$ is a countable (at most) set. The sets $G_l(z)=\{(x,t):\,x\in
R_z,0 \le t \le l,b^t(x)\in M\}$ are finite for all $l\in\R_+$,
and $G(z)=\cup G_l(z)$.

%\medskip

We define the {\em weighted counting function} by
\begin{equation}            \label{coun_cont_eq}
g(z;w)=\sum_{(x,t)\in G(z)}w(x,t).
\end{equation}
The {\em pure counting function} $g_l(z)$ corresponds to the
weight $w=\chi_l$. We have
\begin{equation}      \label{pur_coun_cont_eq}
g_l(z)=|G_l(z)|.
\end{equation}

\medskip

\begin{prop}   \label{crof_cont_prop}
Let $dz,dm$ be finite, lebesgue-class  measures on $Z,M$
respectively; let $dt$ be the lebesgue measure on $\R$. Then there
exist a positive function $r(\cdot)$ on $M\times\R_+$,  determined
by the data 1) and 2), and such that
\begin{equation}      \label{crof_cont_eq}
\int_{Z}g(z;w)dz = \int_{M\times\R_+}w(b^{-t}\cdot m,t)r(m,t)dm\,dt
\end{equation}
\begin{proof}
We define the mapping $f:M\times\R_+\to Z$ by $f=q\circ b^{-t}$.
By conditions 1 and 2, $f$ has full rank almost everywhere. The
pull-back by $f$ of $dz$ is absolutely continuous with respect to
$dmdt$, hence $ f^*(dz)=r(m,t)dm\,dt.$

For $0<l$ set  $w_l(x,t)=w(x,t)\chi_l(t)$, and let $g_l(z;w)$ be
the corresponding counting function. Set
$I_l(w)=\int_{Z}g_l(z;w)dz$. A point, $x\in\Psi$, contributes to
$I_l(w)$ iff $x\in \vp(M\times[0,l])$. Under the change of
variables $dz=d(q\circ\vp(m,t))=r(m,t)dm\,dt$, we have
%
%\begin{equation}      \label{crof_cont_fin_eq}
$$
I_l(w)=\int_{M\times[0,l]}w(b^{-t}\cdot m,t)r(m,t)dm\,dt.
$$
%\end{equation}
%
In the limit $l\to\infty$, we obtain the claim.
\end{proof}
\end{prop}

%\newpage

\subsection{Special cases}  \label{count_spec_sub}
We will discuss a few special cases of
Proposition~\ref{crof_disc_prop} and
Proposition~\ref{crof_cont_prop}. First, the discrete version. The
function $g_n(\tht)$ counts the number of visits in $Y$ of points
$x\in R_{\tht}$ during the first $n$ steps of their journey. Set
$\rho_k=\int_Yr_k(y)dy$, and $R_n=\sum_{k=0}^{n-1}\rho_k$. Then
$\rho_k$ is the volume of $Y_k=T^{-k}(Y)$ with respect to the
measure $\eta^*(d\tht)$. Proposition~\ref{crof_disc_prop} yields
\begin{equation}            \label{pure_disc_coun_eq}
\int_{\Tht}g_n(\tht)d\tht=R_n.
\end{equation}
%%

%\medskip
%\bigskip

In the continuous case the function $g_l(z)$ counts the number of
visits in $M$ of orbits $b^t\cdot x,\,x\in R_{z},$ during the
period $0\le t \le l$. Let $R(l)$ be the volume of the manifold
$N_l\subset \Psi$ with respect to the measure $q^*(dz)$.
Proposition~\ref{crof_cont_prop}  yields
\begin{equation}            \label{pure_cont_coun_eq}
\int_{Z}g_l(z)dz\ =\ R(l).
\end{equation}
%
%\newpage

\section{Bounds on counting functions}  \label{bound}
In this section we analyze the setting of section~\ref{count} from
the measure theoretic viewpoint. This allows us to obtain
pointwise upper bounds on counting functions in a broad spectrum
of situations.
%We will apply these bounds to the billiard setting in section~\ref{exa}.
%

\medskip

Let $X,\mu$ be a finite measure space. Let $f(x;t)$ (for $t\in\R_+$) be a family of nonnegative
$L^1$ functions on $X$. Set
\begin{equation} \label{eqn:croft}
F(t)=\int_X f(x;t) d\mu(x).
\end{equation}
\begin{lem} \label{lem:1cont}
For almost every $x\in X$ there exists $C=C(x)>0$ such that for
arbitrarily large $n\in\N$ there is  $t\geq n$ satisfying $f(x,t)
<  C F(t)$.
\begin{proof}
For $0<C$ and $n\in\N$ let
$$
B_n(C)=\{x\in X:\,  C F(t) < f(x;t)\ \forall t>n \},
$$
and set
$$
B(C)=\bigcup_{n\in\N} B_n(C).
$$
Integrating the inequality above, we obtain $\mu(B_n(C))\leq
C^{-1}$ for any $n$. Thus $\mu(B(C))\leq  C^{-1}$, and hence
$\mu(\cap_{C\in\R_+}B(C))=0$. But $\cap_{C\in\R_+}B(C)\subset X$
is the complement of the set of points $x\in X$ satisfying the
hypothesis of the lemma.
\end{proof}
\end{lem}

\medskip

Let the setting be as in Lemma~\ref{lem:1cont}. In addition, we
suppose that\\
i) the functions $f(x;t)$ are nondecreasing in $t$ and ii)
$F(t)\to\infty$.
\begin{lem} \label{lem:2cont}
Let $\ep >0$ be arbitrary. Then
for almost every $x\in X$ there exists $T=T(x,\ep)>0$ such that for all
$t>T$ we have
\begin{equation}    \label{lem2_eq}
f(x;t)\leq F(t)(1+\log(1+F(t)))^{1+\ep}.
\end{equation}
\begin{proof}
Denote by $f(x;t^-)$ (resp. $F(t^-)$) the limits of
$f(x;s)$ (resp. $F(s)$), as $s\to t$ from the left.
For $n\in\N$ set $t_n = \inf\{t:\, F(t)\geq 2^n\}$.  Then
$F({t_{(n+1)}}^-)\leq 2 F(t_n)$. Let $A_n\subset X$ be the set of points
satisfying the inequality

%It suffices to prove that for almost all $x$ there exists $n_0=n_0(x)$ such that
%for  $n_0<n$ we have
%
\begin{equation}    \label{eqn:lem2}
f(x;t_n^-) \leq \frac 1 2 F({t_n}^-) (1+ \log(1+ \frac12 F({t_n}^-)))^{1+\epsilon}.
\end{equation}
It suffices to prove that the set $\bigcup_{n\in\N} \bigcap_{k>n} A_k$ has full measure.
Indeed, for $x\in A_n$ and $t\in[t_n,t_{n+1})$ we have
$$
f(x;t)\leq f(x;t_{n+1}^-) \leq\frac12 F({t_{n+1}}^-)(1+\log(1+ \frac12 F({t_{n+1}}^-))^{1+\ep}
$$
$$
\leq F(t_n) (1+\log (1+F(t_n)))^{1+\ep} \leq F(t) (1+\log (1+F(t)))^{1+\ep}.
$$
Thus, the points $x\in\bigcup_{n\in\N} \bigcap_{k>n} A_k$ have the
property equation~\eqref{lem2_eq}.

If $B_n\subset X$ is any sequence of sets, we set
$\limsup_{n\to\infty} B_n=\bigcap_{i\in\N} \bigcup_{j>i} B_j$. Let $B_n$ be the complement
of $A_n$ in $X$. Then $\limsup_{n\to\infty} B_n$ is the complement of  $\bigcup_{n\in\N} \bigcap_{k>n} A_k$.
It remains to prove that $\mu(\limsup_{n\to\infty} B_n)=0$.

By Tchebysheff inequality, we have
\begin{equation}  \label{eqn:lem22}
\mu (B_n) \leq 2 (1+\log (1+ \frac12 F({t_n}^-)))^{-(1+\ep)}.
\end{equation}
Set $\mu_n=\mu (B_n)$. Suppose first that $F$ is a continuous function.
Then $F({t_n}^-)=F({t_n})=2^n$. By  equation~\eqref{eqn:lem22}
$$
\mu_n\leq 2 (1+\log (1+ 2^{n}))^{-(1+\ep)},
$$
hence the series $\sum\mu_n$ converges. Since
$$
\mu(\limsup_{n\to\infty} B_n)\le \sum_{n_0}^{\infty}\mu_n
$$
for any $n_0\in\N$, the claim follows.

In general, $F$ need not be continuous. It is thus possible that $t_n=t_{n+1}$
for some $n\in\N$, implying  $B_n=B_{n+1}$. From the series $\sum\mu_n$ we drop
the terms $\mu_n$ such that $B_n=B_{n-1}$. By equation~\eqref{eqn:lem22}, the remaining
terms satisfy
$$
\mu_n \leq 2 (1+\log (1+ 2^{n-2}))^{-(1+\ep)}.
$$
Now the preceding argument applies.
\end{proof}
\end{lem}

%Note that in discrete case all $t_n$ are integers, as they should.

In sections~\ref{exa},~\ref{appli} we will apply these results in
the billiard setting. In section~\ref{exa} we will estimate the
integrals equation~\eqref{eqn:croft}, hence the bounds provided by
Lemmas~\ref{lem:1cont},~\ref{lem:2cont} will be more specific. The
propositions below anticipate these applications.

\begin{prop} \label{appli_prop1}
Let the setting and the  assumptions be as in Lemma \ref{lem:2cont}.
Let $0<\ep$ be arbitrary.

\noindent 1. Let $F(t)=O(t^p)$ for $0<p$.  Then
for almost every $x\in X$ we have $f(x;t)=O(t^{p+\ep})$.

\noindent 2. Let $F(t)=O(e^{at})$  for $0<a$. Then
for almost every $x\in X$ we have $f(x;t)=O(e^{(a+\ep)t})$.
\begin{proof}
The first claim is immediate from Lemma~\ref{lem:2cont} and $(\log t)^{1+\ep}=o(t^{\ep})$.
The second claim follows the same way from $t^{1+\ep}=o(e^{\ep t})$.
\end{proof}
\end{prop}

\medskip

%The variable $t\in\R_+$ in the preceding material will be interpreted as time in
%applications to the billiard flow.
For applications  to the billiard map we need a counterpart of
Proposition~\ref{appli_prop1} for integer-valued time. We state it
below. Its proof is analogous to the proof of
Proposition~\ref{appli_prop1}. Moreover, the discrete time case
may be formally reduced to the continuous time case. We leave
details to the reader.

Let $X,\mu$ be a finite measure space. Let $f(x;n),\,n\in\N$ be a sequence of nonnegative
$L^1$ functions on $X$ such that for every $x\in X$ the numerical sequence $f(x;n)$
is nondecreasing. Set $F(n)=\int_X f(x;n) d\mu$.
\begin{prop} \label{appli_prop2}
Let $0<\ep$ be arbitrary. Then the following claims hold.

\noindent 1. Let $F(n)=O(n^p)$ for $0<p$.  Then
for almost every $x\in X$ we have $f(x;n)=O(n^{p+\ep})$.

\noindent 2. Let $F(n)=O(e^{an})$  for $0<a$. Then
for almost every $x\in X$ we have $f(x;n)=O(e^{(a+\ep)n})$.
\end{prop}
\begin{rem}  \label{appli_rem}
{\em All of the bounds $f(\cdot)=O(\cdot)$ in preceding
propositions are equivalent to the formally stronger bounds
$f(\cdot)=o(\cdot)$.}
\end{rem}

\section{Counting functions for polygonal billiard}   \label{exa}
We will now apply the preceding material to the billiard dynamics.
Our billiard table will be a geodesic polygon either in the
euclidean plane $\RR$, or the round sphere $\unisph$, or the
hyperbolic plane $\hyp$. We refer to \cite{Gut86}, \cite{GT04},
and section~\ref{cover} for the background.

\subsection{Direction counting functions for billiard maps in euclidean polygons}  \label{dir_count_sub}
Let $P\subset\RR$ be a euclidean polygon, and let $T:X(P)\to X(P)$
be the {\em billiard map}. Elements of the phase space $X=X(P)$
are oriented geodesic segments  in $\RR$ with endpoints in $\bo
P$. A segment $x\in X$ ending in a corner of $P$ is {\em
singular}; the element $Tx$ is not well defined. A billiard orbit
$x,Tx,\dots,T^{k-1}x$ is a {\em singular orbit of length $k$} if
$T^{k-1}x$ is the first singular element in the sequence.

\begin{figure}[htbp]
\begin{center}
\input{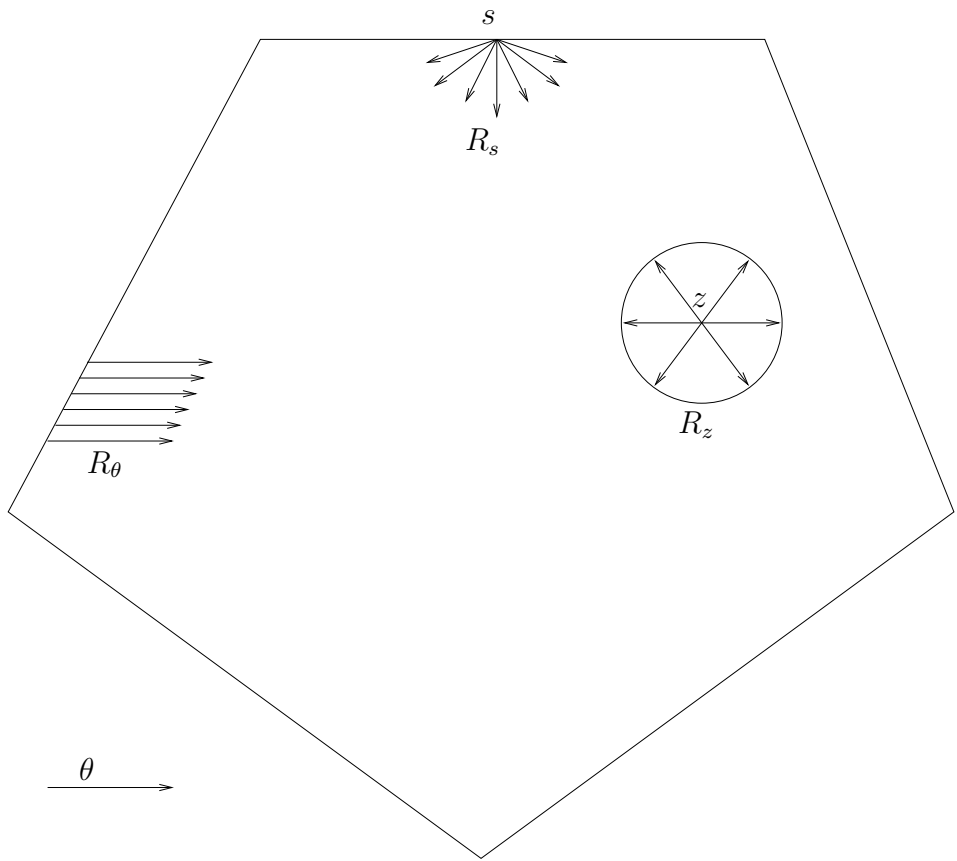}
\caption{Base sets for billiard counting functions}
\label{base_fig}
\end{center}
\end{figure}

Assigning to $x\in X$ its {\em direction}, $\eta(x)\in S^1$, we
obtain a fibration $\eta:X\to S^1$ with fibers $R_{\tht}\subset
X$. See figure~\ref{base_fig}.
%Each $R_{\tht}$ is a convex graph,
%according to the convex structure on $X$ \cite{GT04}.
We define the {\em counting function $gd_{\tht}(n)$ for singular
orbits in direction $\tht$} as the number of phase points $x\in
R_{\tht}$ that yield singular orbits of length $k\le n$.

\begin{thm}   \label{bil_direct_thm}
Let $P\subset\RR$ be  an arbitrary polygon. Let $K(P)$ be the set
of its corners. Let $\al(v)$ be the angle of $v\in K(P)$. Let
$d\tht$ be the lebesgue measure on $S^1$.

Let $K\subset K(P)$. Then
\begin{equation}   \label{bil_direct_eq}
\int_{S^1}\sum_{v\in K}gd_{\tht}(n;v)d\tht = \left(\sum_{v\in
K}\al(v)\right)n.
\end{equation}
\begin{proof}
It suffices to prove the claim for a singleton, $K=\{v\}$. Let
$Y=Y(v)\subset X$ be the set of segments $x\in X$ ending at $v$.
Let $dy$ be the angular measure on $Y$. These data fit into the
setting of section~\ref{count_disc_sub}, and $gd_{\tht}(n;v)$ is
the pure counting function.

Let $B(z,\al)$ be a conical beam of light with apex angle $\al$
emanating from $z\in\RR$. After reflecting in $\bo P$, it splits
into a finite number of beams $B(z_i,\al_i)$ satisfying
$\sum\al_i=\al$. The {\em preservation of  light volume} is due to
the flatness of $\bo P$.

By preceding remark, the functions $r_k(\cdot)$ of
Proposition~\ref{crof_disc_prop} satisfy $r_k(\cdot)\equiv 1$. The
claim now follows from the special case of
Proposition~\ref{crof_disc_prop} considered in
section~\ref{count_spec_sub}.
\end{proof}
\end{thm}

Let $p,q$ be the numbers of corners, obstacles in $P$. Let
$\ka(P)=p+2q-2$. Thus, $P$ is simply connected iff $q=0$ iff
$\ka(P)=p-2$.

\begin{corol}                   \label{bil_direct_cor}
Let $P\subset\RR$ be  an arbitrary polygon. Then
\begin{equation}   \label{bil_direct_eq1}
\int_{S^1}gd_\theta(n)d\theta = \pi\ka(P)n.
\end{equation}
\begin{proof}
Follows from Theorem~\ref{bil_direct_thm} via $\sum_{v\in
K(P)}\al(v)=(p+2q-2)\pi$.
\end{proof}
\end{corol}
%\newpage

%

\subsection{Position counting functions for billiard flows in euclidean polygons}
\label{pos_count_flow}
%
%This section concerns the billiard flow.
Let $P\subset\RR$ be  a  polygon, and let $b^t:\Psi\to\Psi$ be the
billiard flow. See section~\ref{cover} for details. For $z\in P$
and $v\in K(P)$ let $gc_z(l;v)$ be the number of billiard flow
orbits that start from $z\in P$ and wind up at $v$ by time $l$.
Then $gc_z(l)=\sum_{v\in K(P)}gc_z(l;v)$ is the number of singular
billiard orbits of length at most $l$ starting from $z$. This is
the {\em position counting function for the billiard flow} in $P$.

\begin{thm}  \label{pos_count_flow_thm}
Let $P\subset\RR$ be  a euclidean polygon, and let $dz$ be the
lebesgue measure on $P$. Then for any $K\subset K(P)$ we have
\begin{equation}  \label{av_eq3}
2\int_P\sum_{v\in K}gc_z(l;v)dz= \left(\sum_{v\in
K}\al(v)\right)l^2.
\end{equation}
\begin{proof}
It suffices to prove the claim for $K=\{v\}$. We view elements of
$\Psi$ as pairs $z,\tht$ where $z\in P$ is the basepoint, and
$\tht$ is the direction. Let $M=\{(v,\tht):\,(v,-\tht)\in\Psi\}$.
Let $q:\Psi\to P$ be the obvious projection. Its fibers $R_z$ are
the base sets for the counting functions $gc_z(l;v)$. See
figure~\ref{base_fig}. Set $w=\chi_l$. These data satisfy the
assumptions of Proposition~\ref{crof_cont_prop}, and $gc_{z}(l;v)$
is the pure counting function.

We set $dm$ to be the angular measure, and compute the function
$r(m,t)$ in equation~\eqref{crof_cont_eq}. By
Corollary~\ref{pull_back_cor} in section~\ref{cover}, $r=t\chi_l$.
Proposition~\ref{crof_cont_prop} implies the claim.
\end{proof}
\end{thm}

When $K=K(P)$, the left hand side in equation~\eqref{av_eq3} is
the average of the position counting function. The argument of
Corollary~\ref{bil_direct_cor} yields the following.
\begin{corol}                   \label{pos_count_flow_cor}
Let $P\subset\RR$ be  an arbitrary polygon. Then
\begin{equation}   \label{ave_eq4}
2\int_P gc_z(l)dz = \pi\ka(P)l^2.
\end{equation}
\end{corol}

\medskip

%The following estimate is immediate from Corollary~\ref{quad_log_cor}.

%\begin{corol}   \label{power_bound_cor}
%Let $P$ be an arbitrary polygon, and  let $g(z, l)$, $z\in P,\,0<l$ be the counting function
%of billiard orbits joining $z$ with corners.
%Let $\delta>0$ be arbitrary.  Then  for lebesgue almost all $z\in P$ there exists $l(\delta)\in\R$ such that
%for $l>l(\delta)$ we have $g(z, l)<l^{2+\delta}$.
%\end{corol}

%\vspace{10mm}

%\newpage

\subsection{Position counting functions for billiard maps in euclidean polygons}
\label{pos_count_map}
We will now discuss two billiard map analogs of the preceding
example. Let $P\subset\RR$ be a euclidean polygon, and let
$T:X(P)\to X(P)$ be the billiard map.  The phase space $X=X(P)$
consists of pairs $(s,\alpha)$ where $s$ is the arclentgh
parameter on $\bo P$, and $0<\al<\pi$ is the {\em outgoing angle}.
See \cite{Gut86,GT04} and section~\ref{cover} for details. An
orbit $x,Tx,\dots,T^{k-1}x$ is singular, of (combinatorial) length
$k$ if its last segment ends at a corner of $P$.

Let $s\in\bo P$, $v\in K(P)$. Define $GD_s(n;v)$ to be the set of
phase points $(s,\al)\in X$ whose orbits of length less than or
equal to $n$ end at $v$. Set
$$
gd_s(n;v)=|GD_s(n;v)|,\ god_s(n;v)=\sum_{(s,\al)\in
GD_s(n;v)}\sin\al.
$$
The expressions
$$
gd_s(n)=\sum_{v\in K(P)}gd_s(n;v),\ god_s(n)=\sum_{v\in
K(P)}god_s(n;v)
$$
are the {\em pure position counting function} and the {\em optical
position counting function} for the billiard map in $P$.

Let $z\in\RR$ and let $\ga\subset\RR$ be an oriented piecewise
$C^1$ curve. Denote by $d_zs$ the projection of the arclength form
$ds$ of $\ga$ onto the direction perpendicular to the line from
$z$ to $s\in\ga$. The integral
$\int_{\ga}d_zs=|\op(\ga,z)|\le|\ga|$ is the {\em optical length
of $\ga$ viewed from $z$}.

\medskip

Let $z\in P$. Unfolding $k$-segment billiard orbits emanating from
$z$, we obtain a set of linear segments in $\RR$.
 Let $\bo_z(P;k)\subset\RR$ be the curve
traced by their endpoints. We say that $\bo_z(P;k)\subset\RR$ is
the {\em outer boundary of $P$, as viewed from $z$, after $k$
iterates}.
%See section~\ref{cover}.

\begin{thm}  \label{pos_coun_map_thm}
Let $P$ be a euclidean polygon, and let $K\subset K(P)$ be a set
of corners. Then
\begin{equation}  \label{pos_count_map_eq}
\int_{\bo P}\sum_{v\in K}gd_s(n;v)ds=\sum_{v\in
K}\sum_{k=1}^{n}|\bo_v(P;k)|;
\end{equation}

\begin{equation}  \label{pos_count_map_eq1}
\int_{\bo P}\sum_{v\in K}god_s(n;v)ds=\sum_{v\in
K}\sum_{k=1}^{n}|\op(\bo_v(P,k))|.
\end{equation}
\begin{proof}
It suffices to prove the claims for a singleton, $K=\{v\}$. Let
$\eta:X\to\bo P$ be the natural projection. Using the arclength
parametrization, we identify $\bo P$ with the interval $[0,|\bo
P|]\subset\R$. For $0\le s \le |\bo P|$ let
$R_s=\eta^{-1}(s)\subset X$ be the fiber. Then $R_s$ are the base
sets for the counting functions $gd_s(n;v)$, $god_s(n;v)$. See
figure~\ref{base_fig}. Let $Y=Y(v)\subset X$ be the set of phase
points whose $T^{-1}$-orbits emanate from $v$. The assumptions of
section~\ref{count_disc_sub} are satisfied. The weight functions
are $w(s,\al,t)=\chi_n(t)$ and
$w_o(s,\al,t)=\sin\al\cdot\chi_n(t)$ for the two cases at hand.
Let $\vp$ be the angle parameter on $Y$. The measures on $\bo P$
and $Y$ have densities $ds$ and $d\vp$ respectively.

The integrals in the right hand side of
equation~\eqref{crof_disc_eq} are over the curves $\bo_v
(P;k),0\le k\le n-1$. The integrands are $ds(\vp)$ and
$\sin\al\cdot ds(\vp)=d_vs(\vp)$ in respective cases.
\end{proof}
\end{thm}

%\medskip

We will need estimates on lengths and optical lengths.

\begin{lem}  \label{av_pos_map_bound_lem}
For any polygon $P\subset\RR$ there exist $0<c_1 < c_2<\infty$
such that for $n$ sufficiently large
\begin{equation}  \label{pos_bound_map_eq}
c_1n^2 \le \sum_{v\in K}\sum_{k=1}^{n}|\op(\bo_v(P,k))| \le
c_2n^2,\ c_1n^2  \le \sum_{v\in K}\sum_{k=1}^{n}|(\bo_v(P,k))|.
\end{equation}
\begin{proof}
There exist positive constants $d_1,d_2$ and $m_0\in\N$, such that
for any orbit $\ga$ of the billiard map with $m>m_0$ segments, we
have $d_1|\ga| \le m \le d_2|\ga|$ \cite{Gut86}.

Let $v\in\ K(P)$. We will estimate
$\sum_{k=m_0}^{n}|\op(\bo_v(P,k))|$, as $n\to\infty$. Let
$\tht_1\le\tht\le\tht_2$ be the angular parameter for orbits
emanating from $v$; let $r(\tht)$ be the geometric length of the
orbit. Suppose that $r_1\le r(\tht)\le r_2$. Then the optical
length in question is sandwiched between the lengths of circular
arcs of radii $r_1,r_2$ of angular size $\tht_2-\tht_1$. By
preceding remarks, if $k$ is sufficiently large, the bounds
$r_1,r_2$ are proportional to $k$. The total angular size does not
depend on $k$. Hence, for sufficiently large $k$ we have linear
upper and lower bounds on $\sum_{v\in K}|\op(\bo_v(P,k))|$. The
other inequality follows from $|\op(\bo_v(P,k))|\le|\bo_v(P,k)|$.
\end{proof}
\end{lem}

\subsection{Position counting functions for billiard flows in spherical polygons}
\label{sph_pos_count_flow}
%
%This section concerns the billiard flow for a geodesic polygon on
%the round sphere. For concreteness, we will use the unit sphere,
%i. e., $P\subset S^2\subset\R^3$. We will use the analogy with the
%planar case, whuch was discussed in section~\ref{pos_count_flow}.
%The position counting function $gc_z(l)$ is defined as in
%section~\ref{pos_count_flow}, the only difference is that now the
%polygon $P$ is spherical. For a corner $v\in K(P)$ let $gc_z(l;v)$
%be the counting function for orbits that wind up at $v$.
The study is analogous to the planar case discussed in
section~\ref{pos_count_flow}; we will use the same notation
whenever this does not lead to confusion. We denote by $dz$ the
lebesgue measure on $\unisph$, and by $\al(v)$ the angle of a
corner of $P$. Set
\begin{equation}  \label{period_eq1}
\zeta(x) = 1 - \cos x - \frac 2 \pi x.
\end{equation}
\begin{thm}  \label{sph_pos_count_thm}
Let $P\subset\unisph$ be a geodesic polygon, and let $K\subset
K(P)$. Then
\begin{equation}  \label{sph_av_eq3}
\int_P \sum_{v\in K}gc_z(l;v)dz= \left(\sum_{v\in K}\al(v)\right)
\left(\frac 2 \pi l + \zeta(l-\pi \lfloor l/\pi \rfloor)\right).
\end{equation}
\begin{proof}
It suffices to prove the claim when $K=\{v\}$. Let
$M=M(v)\subset\Psi$ be as in section~\ref{pos_count_flow}, and let
$d\al$ be the angular measure on it. The assumptions 1, 3 of
section~\ref{count} are satisfied; the transversality of $b^t\cdot
M$ and $R_z$ may fail for at most a countable set of parameters
$P_{\ex}\subset P$. See Remark~\ref{focus_rem} in
section~\ref{cover}. Hence, condition 2$'$ is fullfilled, and the
results of section~\ref{count_cont_sub} hold. The function
$gc_z(l;v)$ is a pure counting function. The claim now follows
from Proposition~\ref{crof_cont_prop} and
Corollary~\ref{pull_back_cor}.
\end{proof}
\end{thm}

%Let $\vp:A\times\R_+\to P$ be the mapping that assigns to
%$(\al,t)$ the basepoint of $b^t(\al)\in\Psi$. Lifting the orbit
%$b^t(\al)$ to the univeral cover, we identify $\vp$ with the
%exponential map at $v$. See section~\ref{cover}. Hence, $\vp$ is a
%surjective {\bf local?} diffeomorphism on annuli
%$k\pi<t<(k+1)\pi,\,k\in\N$; it collapses the circles $t=k\pi$. The
%pullback of the lebesgue measure on $S^2$ by the exponential map
%is $d\mu=|\sin t|d\al dt$.

%We identify $A\times\R_+$ with the cone $C_v(P)\subset T_vS^2$;
%let $C_v(P;l)$ be the intersection of this cone with the ball of
%radius $l$. Thus  $\mu(C_v(P;l))=\int_{0\le\al\le\al(v)}\int_{0\le
%t \le l}|\sin t| d\al dt$. Evaluating the integral and using
%Proposition~\ref{crof_cont_prop}, we obtain the claim.

Let $\ka(P)$ be as in section~\ref{dir_count_sub}.

\begin{corol}                   \label{sph_pos_count_cor}
Let $P\subset\unisph$ be  an arbitrary polygon. Then
\begin{equation}   \label{sph_av_eq4}
\int_P gc_z(l)dz = \left( \ka(P)\pi + \ar(P) \right) \left(\frac 2
\pi l + \zeta(l-\pi \lfloor l/\pi \rfloor)\right).
\end{equation}
\begin{proof}
For a spherical polygon we have $\sum_{v\in K(P)}\al(v)=\ar(P) +
\ka(P)\pi$. Substitute this into equation~\eqref{sph_av_eq3}.
\end{proof}
\end{corol}
%\newpage
%
\subsection{Position counting functions for  billiard flows in hyperbolic polygons}
\label{hyp_pos_count_flow}
%
%This section concerns the billiard flow for a geodesic polygon in
%the hyperbolic plane. We define the position counting function
%$gc_z(l)$ as in section~\ref{pos_count_flow}, taking into account
%that $P\subset\HH$. For a corner $v\in K(P)$ let $gc_z(l;v)$ be
%the counting function for orbits that wind up at $v$.
Our treatment and our notation are modelled on
section~\ref{sph_pos_count_flow}. We denote by $dz$ the lebesgue
measure on $\HH$, and by $\al(v)$ the angles of corners.
%
%\begin{equation}  \label{hyp_period_eq1}
%\sig(x) = 1 - \cosh x - \frac 2 \pi x.
%\end{equation}
%
\begin{thm}  \label{hyp_pos_count_thm}
Let $P\subset\HH$ be a geodesic polygon, and let $K\subset K(P)$.
Then
\begin{equation}  \label{hyp_av_eq3}
\int_P \sum_{v\in K}gc_z(l;v)dz= \left(\sum_{v\in K}\al(v)\right)
\cosh l.
\end{equation}
\begin{proof}
We repeat verbatim the proof of Theorem~\ref{sph_pos_count_thm},
and use claim $2$ in Corollary~\ref{pull_back_cor}.
\end{proof}
\end{thm}

%It suffices to prove equation~\eqref{hyp_av_eq3} when $K=\{v\}$.
%Let $A\subset\Psi$ be the set of phase points based at $v$, and
%let $d\al$ be the angular lebesgue measure on it. Let
%$\vp:A\times\R_+\to P$ be the mapping that assigns to $(\al,t)$
%the basepoint of $b^t(\al)\in\Psi$. Lifting the orbit $b^t(\al)$
%to the univeral cover, we identify $\vp$ with the exponential map
%at $v$. See section~\ref{cover}. The pullback of the lebesgue
%measure on $\HH$ by the exponential map is $d\mu=\sinh td\al dt$.

Let $\ka(P)$ be as in section~\ref{dir_count_sub}.

\begin{corol}                   \label{hyp_pos_count_cor}
Let $P\subset\HH$ be  a polygon. Then
\begin{equation}   \label{hyp_av_eq4}
\int_P gc_z(l)dz = \left( \ka(P)\pi - \ar(P) \right)\cosh l.
\end{equation}
\begin{proof}
Repeat the argument of Corollary~\ref{sph_pos_count_cor}; use the
formula $\sum_{v\in K(P)}\al(v)=\ka(P)\pi-\ar(P)$ relating the
angles and the area of geodesic polygons in $\hyp$.
\end{proof}
\end{corol}

%For a geodesic polygon in $\HH$ we have $\sum_{v\in
%K(P)}\al(v)=\ar(P) + \ka(P)\pi$. Substitute this into
%equation~\eqref{hyp_av_eq3}.

\section{Relating partial complexities and counting functions}\label{comp_count}
%
%The concept of {\em complexity in dynamics} involves partitioning
%the {\em phase space} into elementary parts (atoms), coding finite
%dynamical orbits by the atoms that they visit, and studying the
%number of these codes, as orbits become longer and longer. The
%complexity of a dynamical system corresponding to a natural
%partition yields valuable information about the dynamics.%
%
%As an example, we consider the billiard in a square, $P$. The four
%sides, $a,b,c,d$ define a partition of the billiard phase. The
%corresponding coding assigns to a billiard orbit the sequence of
%sides  of $P$ that it has visited. Running billiard orbits in $P$
%for a finite time, we generate a finite set of words on the
%alphabet  $a,b,c,d$. The cardinality of this set as a function of
%the length of the time interval is the complexity in question; its
%asymptotics as the time  interval goes to infinity provides useful
%information about the set of billiard orbits in $P$.
%
%The example above is a special case of a natural class of
%examples, as follows. We let $P$ be an arbitrary planar polygon,
%and partition the billiard phase space by atoms defined by the
%sides of $P$. This defines the {\em side complexity} of the
%billiard in $P$. We denote it by $f_P(n),\,n\in\N$. Basic
%questions about $f_P(n)$ are open \cite{Gut03}. For instance,
%$f_P(n)$ is not known when $P$ is a small, but generic,
%perturbation of the square.
%
%\vspace{2mm}
%
In this section we establish a framework that will allow us to
study the complexity of a wide class of dynamical systems. Our
motivation comes from the billiard dynamics. In fact, polygonal
billiard is the target of applications for our results. The
framework is more general, however. The following observations
served as our guiding principles. First, natural partitions of the
billiard-type systems are geared to the singularities. Second, the
billiard dynamics satisfies a certain convexity property that is
instrumental in the study of complexity. These principles are
manifest in the framework of {\em piecewise convex
transformations} \cite{GT04}.

%\medskip

There are two approaches to the billiard dynamics: The {\em
billiard flow} and the {\em billiard map}. See
section~\ref{cover}. The framework of piecewise convex
transformations is geared to the billiard map. We begin by
establishing its counterpart for flows.

\subsection{Piecewise convex transformations and piecewise convex flows}  \label{convex_sub}
A piecewise convex transformation is a triple $(X,\Ga,T)$, where
$X$ is a two-dimensional convex cell complex, $\Ga\subset X$ is
the graph formed by the union of one-cells, and $T:X\to X$ is an
invertible map, regular on the two-cells of the complex, and
compatible with the convex structure \cite{GT04}.
%They are convex, and $T^n$ is regular on these regions.

Let $\Psi$ be a compact manifold, with boundary and corners, in
general. Let $b^t:\Psi\to\Psi$ be a flow, possibly with
singularities; let $X\subset\Psi$ be a cross-section. We will
assume that the singular set of the flow is contained in $X$. For
$z\in X$ let $\tau_+(z),\tau_-(z)$ be the times when $z\in\Psi$
first reaches $X$ under $b^t,b^{-t}$ for $0<t$. We assume that for
any $z\in\Psi\setminus X$ there is $0<\ep=\ep(z)$ such that
$b^t(z)$ is regular for $|t|<\ep$.

A {\em piecewise convex flow} is determined by the following data:
A flow, $b^t:\Psi\to\Psi$, a cross-section, $X\subset\Psi$, and
the structure of a convex cell complex on $X$, compatible with the
poincare map. Billiard flows for polygons on surfaces of constant
curvature are piecewise convex flows \cite{GT04}.

\subsection{Partial complexities for  maps and flows}  \label{complex_sub}
Let $(X,\Ga_n,T^n)$ be the  iterates of a piecewise convex
transformation $(X,\Ga,T)$.\footnote{ They are piecewise convex
transformations as well \cite{GT04}. } Let $F(\Ga_n)$ be the
finite set of open faces of $\Ga_n$; these are the continuity
regions for $T^n$. The function $f(n)=|F(\Ga_n)|$ is the (full)
complexity of $(X,\Ga,T)$.

Let $R\subset X$ be a closed subset. Set
%
%\begin{equation}  \label{part_comp_eq}
$$
F_R(n)=\{A\in F(\Ga_n):A\cap R\ne\emptyset\}.
$$
%\end{equation}
%
%\medskip
%
\begin{defin} \label{part_comp_def}
The function $f_R(n)=|F_R(n)|$ is the {\em partial complexity of
the piecewise convex transformation $(X,\Ga,T)$ based on the
subset $R$.}
\end{defin}

%\medskip

Let $b^t:\Psi\to\Psi$ be a piecewise convex flow, and let
$R\subset\Psi$ be a closed, convex set transversal to the flow.
For $0<l$ let $O_R(l)$ be the set of regular flow orbits of length
$l$ starting from $R$. Let $\al_0,\al_1\in O_R(l)$. A {\em
homotopy} is a continuous family of regular orbits $\al_p\in
O_R(l),\,0\le p \le 1,$ interpolating between $\al_0,\al_1$. We
will say, for brevity, that the orbits $\al_0,\al_1$ are {\em
$R$-homotopic}.
%This yields an equivalence relation on  $O(l)$ (resp. $O_R(l)$);
We denote by $H_R(l)$ the set of $R$-homotopy classes.

%Generally, the combinatorial length is not constant under a homotopy.

%An orbit $x_1,x_2=T(x_1),\dots,x_n$ of the billiard map determines
%a billiard flow orbit $\al(t)=(z(t),\theta(t)),\,0\le t \le l,$
%comprized of $n$ full segments. We say that
%$l$ (resp. $n$) is the {\em geometric (resp. combinatorial) length} of the orbit.
%The two lengths are equivalent in the following sense \cite{Gut86,Gut96}.
%There are positive constants $d_1,d_2,l_0$ depending on $P$ so that
%for $l_0\le l$ we have
%
%\begin{equation}  \label{comb_geom_eq1}
%d_1\,n \le l \le d_2\,n.
%\end{equation}
%
%If $\al$ is a flow orbit of the geometric length $l$ and with $n$
%complete segments, then for  sufficiently large $l$ we have
%
%\begin{equation}  \label{comb_geom_eq2}
%d_1\,n\le l \le d_2\,n  + 2\,\diam P.
%\end{equation}
%

%{\bf There is a hidden figure here. Maybe we dont need it; decide later on}
%
%\begin{figure}[htbp]
%\begin{center}
%\input{fig1.pstex_t}
%\caption{Geometric length and combinatorial length of a billiard orbit}
%\label{orb_fig}
%\end{center}
%\end{figure}

%Let $R\subset\Psi$ be a set satisfying the requirements of Definition~\ref{part_comp_def},
%and {\em transversal to the billiard flow}.
%For $0<l$ let $A(l)$ (resp. $A_R(l)$) be the set of all (resp.
%starting from $R$) regular orbits of length $l$.

\begin{defin} \label{flow_comp_def}
The function  $h_R(l)=|H_R(l)|$ is the {\em partial complexity
(based on $R$) of the  piecewise convex flow} $b^t:\Psi\to\Psi$.
\end{defin}

%For brevity, we will refer to $f_R(n),h_R(l)$ as the {\em
%$R$-complexity} in the discrete, continuous cases respectively.
%In order to distinguish between the two cases, we will call $f_R(n)$
%(resp. $h_R(l)$) the $R$-complexity of the {\em map} (resp. {\em flow}).
%Our sets $R$ will have geometric meanings
%corresponding to a direction, footpoint, etc. We will then speak of the
%(discrete/continuous) direction complexity,  position complexity, etc.

%
%\medskip

\begin{figure}[htbp]
\begin{center}
\input{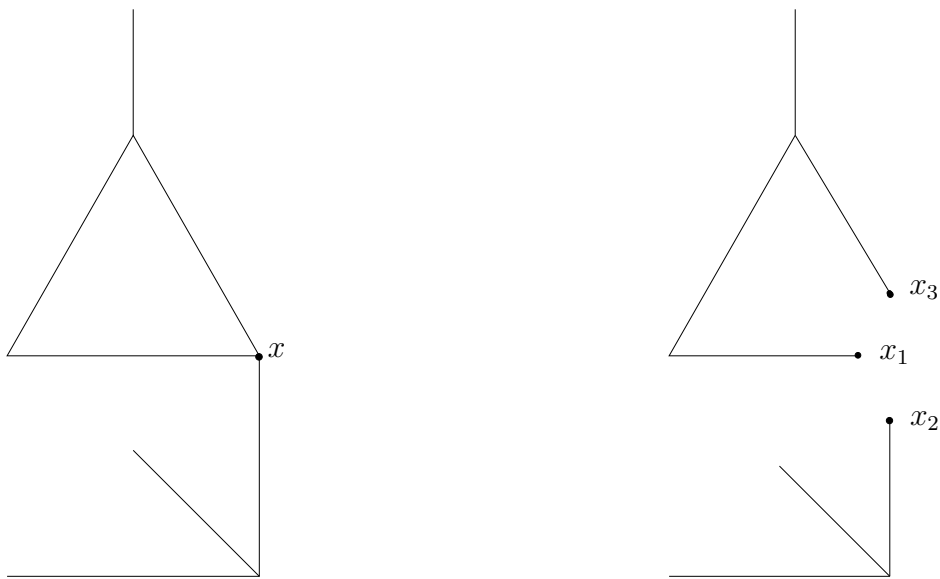}
\caption{Removing a vertex in a graph} \label{orb_fig2}
\end{center}
\end{figure}

\subsection{Partial complexities and counting functions}  \label{count_sub}
In what follows we assume that $R\subset\Psi$ is a {\em convex
graph} without isolated vertices. For $x\in R$ its {\em valence}
$\val(x)$  is  the number of edges of $x$ minus one. In
particular, if $x$ is an interior point of an edge, then
$\val(x)=1$. Set $\val(R)=\max_{x\in R}\val(x)$. We endow
$R\setminus\{x\}$ with the graph structure where $x$ is replaced
by $1+\val(x)$ vertices; each of them is the endpoint of a unique
edge. If $x,y,z,\ldots\in R$ are distinct points, then the
inductively defined graph structure on $R$ without $x,y,z,\ldots$
does not depend on the order of removing these points. We will
denote this graph by $R\setminus\{x,y,z,\ldots\}$. See
figure~\ref{orb_fig2} for an illustration.

Let $E(R)$ and $V(R)$ be the sets of edges and vertices, and let
$c(R)$ be the number of connected components of the graph. Let
$h_i=h_i(R)$ be the betti numbers of $R$, and set
$\chi(R)=|V(R)|-|E(R)|$.  Then $c(R)=h_0,\,\chi(R)=h_0-h_1$.

\begin{lem} \label{count_fun_lem}
Let $R$ be a finite graph, and let $x_1,\dots,x_p\in R$ be
distinct points. Then
\begin{equation}   \label{count_fun_eq}
 \chi(R) + \sum_{i=1}^p \val(x_i) \leq c(R\setminus\{x_1,\ldots,x_p\}) \leq c(R) +
\sum_{i=1}^p \val(x_i).
\end{equation}
If $R$ is a forest, then the bound on the right in
equation~\eqref{count_fun_eq} becomes an equality.
%
%\begin{equation}   \label{count_fun_eq2}
%c(R\setminus\{x_1,\ldots,x_p\}) = c(R) + \sum_{i=1}^p \val(x_i).
%\end{equation}
%
\begin{proof}
It suffices to prove the claims when $R$ is connected, and we
remove a single vertex, $x$. Equation~\eqref{count_fun_eq} becomes
\begin{equation}   \label{count_fun_eq1}
 \chi(R) + \val(x) \leq c(R\setminus\{x\}) \leq c(R) + \val(x).
\end{equation}
We have
$|V(R\setminus\{x\})|=|V(R)|+\val(x),\,|E(R\setminus\{x\})|=|E(R)|$,
and $\chi(R\setminus\{x\})=\chi(R)+\val(x)$. Equivalently, we have
$\chi(R\setminus\{x\})=h_0(R)+\val(x)-h_1(R)$ and
$h_0(R\setminus\{x\})=h_0(R)+\val(x)+\left(h_1(R\setminus\{x\})-h_1(R)\right)$.
The former (resp. latter) identity implies the left (resp. right)
inequality in equation~\eqref{count_fun_eq1}.

When $R$ is a tree, we have $c(R\setminus\{x\}) = c(R) + \val(x)$,
and the remaining claim follows.
\end{proof}
\end{lem}

%\medskip

We will introduce counting functions for singular orbits of the
billiard map and the billiard flow.

By definition, an orbit $\al=\{b^t(z), 0\le t \le l\},$ does not
pass through singular points in $\Psi$. It is regular if it does
not contain any  singular points in $\Psi$; it is {\em singular}
if one of its endpoints is singular. The set  $S_R(l)$ of singular
orbits of length at most $l$, based in $R$, is finite. The
quantities $gc_R(l)=|S_R(l)|$ and $gd_R(n)=|R\cap\Ga_n|$ are the
{\em counting functions for singular orbits based in $R$} for the
flow and the map respectively.
%{\bf illustrate definition
% withfigure or example}

%
%\begin{defin}   \label{count_fun_def}
%Let $R\subset\Psi$ be a transversal graph, and let $0\le l$.
%The function  $gc_R(l)=|S_R(l)|$ (resp. $gd_R(n)=|R\cap\Ga_n|$)
%is the {\em counting function of singular flow (resp. map) orbits  based in $R$}.
%\end{defin}
%
%For brevity, we will speak of the  map (resp. flow) counting functions based on $R$.
%

Now we will relate partial complexities and counting functions. We
do this for a piecewise convex flow $b^t:\Psi\to\Psi$ and for a
piecewise convex transformation $(X,\Ga,T)$. In both cases the
partial complexity is based on a $1$-dimensional subset, say $R$.
Recall that $gc_R(l),gd_R(n)$ are the respective counting
functions, and $h_R(l),f_R(n)$ are the respective complexities. We
will refer to these situations as {\em the continuous  case and
the discrete case} respectively.

\begin{prop}  \label{comp_count_prop}
Let the setting be as above. Then the following statements hold.

\noindent 1. In the continuous  case there exist $h_0\in\N$ and
$l_0\in\R_+$ such that $h_R(l) = h_0 + gc_R(l)$ for $l_0\le l$. 2.
In the discrete  case there exist $f_0,n_0\in\N$ such that for
$n_0\le n$ we have $f_R(n) = f_0 + gd_R(n)$.
\begin{proof}
In both cases the graph $R$ is equipped with a tower of finite
sets, say $X(l)$ and $X_n$ respectively. Let $X_\infty\subset R$
be their union. We will compare the number of connected components
of graphs $R\setminus X(l),R\setminus X_n$ with the cardinalities
of these sets.

We consider the discrete case, leaving the continuous case to the
reader. Let $m<n$ be any pair of natural numbers. By (the proof
of) Lemma~\ref{count_fun_lem},
%
%\begin{equation}   \label{compare_eq}
$$
c(R\setminus X_n)-c(R\setminus X_m)=\left[h_1(R\setminus X_n)-
h_1(R\setminus X_m)\right]+\sum_{x\in X_n\setminus X_m}\val(x).
$$
%\end{equation}
%

We have $h_1(R\setminus X_n)\le h_1(R\setminus X_m)$; the
inequality holds iff $X_n\setminus X_m$ breaks cycles in
$R\setminus X_m$. Since the sequence $h_1(R\setminus X_k)\in\N$ is
nonincreasing, it stabilizes. Thus, there exists $n_1\in\N$ such
that for $n_1\le m<n$ we have $h_1(R\setminus X_n)=h_1(R\setminus
X_m)$.

The set of points $x\in R$ satisfying $1< \val(x)$ is finite.
Thus, there exists $n_2\in\N$ such that if $n_2\le k$ and $x\in
X_\infty\setminus X_k$, then $\val(x)=1$.

Set $n_0=\max(n_1,n_2)$. Then for $n_0\le m<n$ the above equation
yields $c(R\setminus X_n)-c(R\setminus X_m)=|X_n\setminus X_m|.$
Specializing to $m=n_0$, we obtain
$f_R(n)=\left(f_R(n_0)-gd_R(n_0)\right)+gd_R(n)$.
\end{proof}
\end{prop}

%If $\al\in S_R(l)$ (resp. $x\in R\cap\Ga_n$), then $\val(\al(0))$
%(resp. $\val(x)$) is the {\em valence of the corresponding
%singular orbit}.

\section{Bounds on partial complexities for the billiard}  \label{appli}
We will use the preceding material to derive bounds on partial
complexities for the polygonal billiard.

\subsection{Direction complexities for billiard maps in euclidean polygons}  \label{dir_comp_sub}
We use the setting and the notation of
section~\ref{dir_count_sub}. For a  polygon  $P$  and a direction
$\theta$, we denote by $fd_{\tht}(n)$ the partial complexity with
base $R_{\tht}$. This is the {\em complexity in direction $\tht$}.

\begin{corol}       \label{dir_est_cor1}
For lebesgue almost all directions $\theta$ there is $C=C(\theta)$
and there are arbitrarily large $n$ such that $fd_\theta(n) \leq
Cn$.
\begin{proof}
Each $R_{\tht}$ is a convex graph in the phase space \cite{GT04}.
By Lemma~\ref{lem:1cont} and Corollary~\ref{bil_direct_cor}, the
counting functions $gd_\theta(n)$ have the desired properties. By
the second claim of Proposition~\ref{comp_count_prop}, the
directional complexities do as well.
\end{proof}
\end{corol}
\begin{corol}          \label{dir_est_cor2}
For any $\ep>0$ and almost every direction $\theta$ we have
$fd_\theta(n) = O(n^{1+\ep})$.
\begin{proof}
The proof goes along the lines of the proof of
Corollary~\ref{dir_est_cor1}. Instead of Lemma~\ref{lem:1cont}, we
use Proposition~\ref{appli_prop2} (the first claim).
\end{proof}
\end{corol}

\subsection{Position complexities for billiard flows in euclidean polygons}
\label{pos_comp_flow}
Let $P$ be a euclidean polygon, and let $z\in P$ be any point. We
consider the billiard flow in $P$, and use the setting of
section~\ref{pos_count_flow}. Thus, $gc_z(l)$ is the position
counting function for orbits emanating from $z$. We denote by
$h_z(l)$ the corresponding partial complexity.
\begin{corol}       \label{quad_est_cor}
For almost every point $z$ there is a positive number $C=C(z)$
such that  $h_z(l) \leq C l^2$ for arbitrarily large $l$.
\begin{proof}
The sets $R_z$ satisfy the assumptions of
section~\ref{comp_count}. The claim follows from
Lemma~\ref{lem:1cont}, Corollary~\ref{pos_count_flow_cor} and the
continuous case in Proposition~\ref{comp_count_prop}.
\end{proof}
\end{corol}

\begin{corol}          \label{quad_log_cor}
For any $\ep>0$ and almost every $z\in P$ we have $h_z(l) =
O(l^{2+\ep})$.
\begin{proof}
The proof is similar to the preceding argument, and we use the
first claim in Proposition~\ref{appli_prop1} instead of
Lemma~\ref{lem:1cont}.
\end{proof}
\end{corol}
\subsection{Position complexities for billiard maps in euclidean polygons}
\label{pos_comp_map}
This is the billiard map analog of the preceding example. Let $P$
be a euclidean polygon, and let $s\in\bo P$. We use the setting of
section~\ref{pos_count_map}. There we have defined the counting
functions $gd_s(n),god_s(n)$. Let $f_s(n)$ be the partial
complexity corresponding to $gd_s(n)$. This is the {\em position
complexity for the billiard map}.
\begin{corol}  \label{pos_coun_map_cor}
%Let $P\subset\RR$ be any polygon. Then for almost all $s\in\bo P$
%and any $0<\ep$ we have $f_s(n)=O(n^{2+\ep})$.
Let $P\subset\RR$ be a polygon such that $\sum_{v\in
K}\sum_{k=1}^{n}|\bo_v(P;k)|$ has a quadratic upper
bound.\footnote{This holds if $P$ is a rational polygon
\cite{M90}.} Then for almost all $s\in\bo P$ we have
$f_s(n)=O(n^{2+\ep})$ for any $0<\ep$.
\begin{proof}
The sets $R_s\subset X$ satisfy the assumptions of
section~\ref{comp_count}. We use Theorem~\ref{pos_coun_map_thm},
Lemma~\ref{av_pos_map_bound_lem}, and  apply
Proposition~\ref{comp_count_prop}.
\end{proof}
\end{corol}
The estimate of Corollary~\ref{pos_coun_map_cor} on  $f_s(n)$ is
conditional, because in general we have no efficient upper bound
on $\sum_{k=1}^{n}|\bo_v(P;k)|$.

\subsection{Position complexities for billiard flows in spherical polygons}
\label{sph_pos_comp}
We use the setting of section~\ref{sph_pos_count_flow}. For a
spherical polygon, $P\subset\sph$,  and $z\in P$, let $h_z(l)$ be
the position complexity.

\begin{corol}       \label{sph_est_cor1}
For almost every point $z\in P$ there is $C=C(z)$ and there are
arbitrarily large $l$ such that $h_z(l) \leq C l$.
\begin{proof}
The sets $R_z$ satisfy the assumptions of
section~\ref{comp_count}. We use Lemma~\ref{lem:1cont},
Corollary~\ref{sph_pos_count_cor},  and
Proposition~\ref{comp_count_prop}.
\end{proof}
\end{corol}

\begin{corol}          \label{sph_est_cor2}
For any $\ep>0$ and almost every $z\in P$ we have $h_z(l) =
O(l^{1+\ep})$.
\begin{proof}
See the proof of Corollary~\ref{quad_log_cor}.
\end{proof}
\end{corol}

\subsection{Position complexities for billiard flows in hyperbolic  polygons}
\label{hyp_pos_comp}
This material is the hyperbolic plane counterpart of
section~\ref{pos_count_flow}, and we use the setting of
section~\ref{hyp_pos_count_flow}.
\begin{corol}       \label{hyp_est_cor}
Let $P\subset\HH$ be a geodesic polygon, let $z\in P$,  and let
$h_z(l)$ be the position complexity. Then for almost every point
$z\in P$ we have $h_z(l)=O(e^{(1+\ep)l})$.
\begin{proof}
We verify that the sets $R_z$ satisfy the assumptions of
section~\ref{comp_count}, and mimick the proof of
Corollary~\ref{quad_log_cor}; we use
Corollary~\ref{hyp_pos_count_cor}, Proposition~\ref{appli_prop1},
and the continuous case of Proposition~\ref{comp_count_prop}.
\end{proof}
\end{corol}

%\newpage

%\newpage

\section{Appendix: Covering spaces for polygonal billiards}   \label{cover}
Let $M$ be a simply connected surface of constant curvature, and
let $P\subset M$ be a connected geodesic polygon. We normalize the
metric so that the curvature is either zero ($M=\RR$), or one
($M=\unisph$), or minus one ($M=\HH$).
%A spherical polygon
%$P\subset\unisph$ is {\em admissible} if it is contained in an
%open hemisphere \cite{GT04}. If $M=\RR$ or $M=\HH$, then any
%geodesic polygon is admissible. \marginpar{only because
%\cite{GT04} had the assumption} We will assume from now on that
%our polygon is admissible.

Let $A$ be the set of sides in $P$. We will denote its elements by
$a,b,\dots$. For a side, say $a\in A$, let $s_a\in\iso(M)$ be the
corresponding geodesic reflection.
%Let $\iso_P(M)\subset\iso(M)$ be the group of
%isometries generated by $\{s_a:a\in A\}$.
We associate with $P$ a {\em Coxeter system} $(G,A)$ \cite{Dav83}.
We denote by $\sig_a,\sig_b,\dots\in G$ the elements corresponding
to $a,b,\ldots\in A$. They generate $G$. The defining relations
are $\sig_a^2=1$ and $(\sig_a\sig_b)^{n(a,b)}=1$; the latter arise
only for the sides $a,b$ with a common corner if the angle,
$\theta(a,b)$, between them is $\pi$-rational. In this case
$n(a,b)$ is the denominator of $\theta(a,b)/\pi$. Otherwise
$n(a,b)=\infty$.

To any ``generalized polyhedron" $P$ corresponds a topological
space $\cc$ endowed with several structures, and a Coxeter system
\cite{Dav83}. Our situation fits into the framework of
\cite{Dav83}, and we apply its results. First, $\cc$ is a
differentiable surface. Second, $\cc$ is tiled by subsets
$P_g,g\in G$, labelled by elements of the Coxeter group $G$; we
call them the {\em tiles}, and identify $P_e$ with $P$. The group
$G$ acts on $\cc$ properly discontinuously, preserving the tiling:
$g\cdot P_h=P_{gh}$.

Since $P_e$ is identified with $P\subset M$, it inherits from $M$
a riemannian structure. The action of $G$ is compatible with this
structure, and extends it to all of $\cc$. This riemannian
structure generally has cone singularities at vertices of the
tiling $\cc=\cup_{g\in G}P_g$.\footnote{Each vertex, $v$,
corresponds to a corner of $P$. The metric at $v$ is regular iff
the corner angle is $\pi/n,\,n=2,3,\dots$.} Around other points
this riemannian structure is isometric to that of $M$; in
particular, except for cone points, $\cc$ has constant curvature.
The group $G$ acts on $\cc$ by isometries.
\begin{defin}  \label{uni_cov_def}
The space $\cc$ endowed with the riemannian structure, the
isometric action of $G$ and the $G$-invariant tiling
$\cc=\cup_{g\in G}P_g$ is the {\it universal covering space} of
the geodesic polygon $P\subset M$.
\end{defin}

If $X$ is a riemannian manifold (with boundary and singularities,
in general), we denote by $TX=\cup_{x\in X}T_xX$ its {\em unit
tangent bundle}.
%We have $TX=\cup_{x\in X}T_xX$, where $T_xX$
%consists of vectors with base-point $x$.
The classical construct of {\em geodesic flow}, $G^t_X:TX\to TX$,
%which is a classical construct if $X$ has no singularities,
extends to manifolds with boundaries and singularities. In
particular, $G^t_X$ makes sense when $X=M,P$, or $\cc$. Another
classical construct, the {\em exponential map}, also extends to
our situation. For $x\in X$ as above, and $(v,t)\in
T_xX\times\R_+$, we set $\exp_X(v,t)\in X$ be the base-point of
$G^t_X(v)$. We will use the notation $\exp_X^x$ to indicate that
we are exponentiating from the point $x$. If $X$ is nonsingular,
then $\exp_X^x:T_xX\times\R_+\to X$ is a differentiable mapping.
For $X$ with singularities, such as our $P$ and $\cc$, the maps
$\exp_X^x$ are defined on proper subsets of $T_xX\times\R_+$;
these subsets have full lebesgue measure. Generally, the maps do
not extend by continuity to all of $T_xX\times\R_+$.

Let $X,Y$ be nonsingular riemannian manifolds of the same
dimension; let $\vp:X\to Y$ be a local isometry. It induces a
local diffeomorphism $\Phi:TX\to TY$ commuting with the geodesic
flows: $\Phi\circ G^t_X = G^t_Y\circ\Phi$. The exponential maps
commute as well: $\vp\circ\exp_X^x = \exp_Y^{\vp(x)}\circ\,
d_x\vp$. These relationships hold, in particular, for coverings of
nonsingular riemannian manifolds. Suitably interpreted, they
extend to (branched) coverings of riemannian manifolds with
boundaries, corners, and singularities. In our case $X=\cc$, while
$Y=M$, or $Y=P$. We will now define the mappings $f:\cc\to
P,\,F:T\cc\to TP$  and $\vp:\cc\to M,\,\Phi:T\cc\to TM$.

The identification $P_e=P$ defines $f,\,\vp$ on $P_e$. To extend
them to all of $\cc$, we use the tiling $\cc=\cup_{g\in G}P_g$ and
the actions of $G$ on $\cc$ and $M$. In order to distinguish
between these actions, we will denote them by $g\cdot x$ and
$g(x)$ respectively.
%To simplify notation, we suppress $f,\,\vp$ in
%$f:P_e\to P,\vp:P_e\to M$. Let $z\in P_g,\,g\ne e$.
Then there is a unique $x\in P_e$ such that $z=g\cdot x$. We set
$f(z)=x\in P$ and $\vp(z)=g(x)\in M$. By basic properties of
Coxeter groups \cite{Dav83}, the mappings $f,\vp$ are well
defined. Moreover, $f:\cc\to P$ and $\vp:\cc\to M$ are the unique
$G$-equivariant mappings which are identical on
$P_e$.\footnote{The action of $G$ on $P$ is trivial.} By
construction, both mappings are continuous; they  are
diffeomorphisms in the interior of each tile, $P_g\subset\cc$, and
on the interior of the union of any pair of adjacent tiles.

The potential locus of non-differentiability for both $f$ and
$\vp$ is the set $V$ of vertices in the tiling $\cc=\cup_{g\in
G}P_g$. We have $V=f^{-1}(K(P))$ where $K(P)$ is the set of
corners of $P$. By equivariance, $\vp(V)=\cup_{g\in
G}g(K(P))\subset M$.\footnote{ The representation $M=\cup_{g\in
G}g(P)$ is not a tiling, in general.} There are two kinds of
points in $V$: vertices coming from the corners of $P$ with
$\pi$-rational and $\pi$-irrational angles. Their cone angles are
integer multiples of $2\pi$ and are infinite respectively.
Vertices $v\in V$ with cone angle $2\pi$ are, in fact, regular
points in $\cc$, and the mappings $f,\,\vp$ are both regular
there. Around a vertex $v$ with cone angle $2k\pi>2\pi$ the
mapping $\vp$ is differentiable, but not a diffeomorphism; it is
locally conjugate to $z\mapsto z^k$. Near such a vertex, $\vp$ is
a branched covering of degree $k$. At a vertex with infinite cone
angle, $\vp$ has infinite branching.

\begin{rem}  \label{focus_rem}
{\em The set $\vp(V)\subset M$ is countable. (It is finite iff the
group generated by geodesic reflections in the sides of $P$ is a
finite Coxeter group. Typically, $\vp(V)\subset M$ is a dense,
countable set.) Let $M=\unisph$, and let $z\mapsto z'$ denote the
antipodal map. Set $F=P\cap\left(\vp(V)\cup(\vp(V))'\right)$.
Points  of $F$ are exceptional, in the following sense. Let $z\in
P$ be such that the beam $R_z$ of billiard orbits emanating from
$z$ contains a sub-beam focusing at a corner of $P$. Then $z\in
F$. This follows from Proposition~\ref{straigh_prop} below.

Thus, $F$ contains all points $z\in P$ for which the
transversality assumption in Condition 2 of
section~\ref{count_cont_sub} fails. Since $F$ is countable, the
set of exceptional parameters has measure zero, and Condition 2$'$
is satisfied. See Remark~\ref{except_rem} in section~\ref{count}.}
\end{rem}

Furthermore, the mappings $f$ and $\vp$ are local isometries. They
are isometries on every tile $P_g\subset\cc$; we have $f(P_g)=P$,
$\vp(P_g)=g(P)\subset M$. Let $g\cdot a$ be a side of $P_g$, let
$h=\sig_ag$ and let $P_h$ be the adjacent tile. The maps $f:P_g\to
P,P_h\to P$ and $\vp:P_g\to g(P),P_h\to h(P)$ are coherent around
the common (open) side $g\cdot a$. The map $f$ is never an
isometry on $P_g\cup P_h$; for $\vp$ this is the case iff the
interiors of $g(P),h(P)$ are disjoint in $M$. The latter generally
fails for nonconvex $P$.

By coherence of $f$ and $\vp$ across the sides separating adjacent
tiles, we lift them to the tangent bundles, obtaining the mappings
of unit tangent bundles $F:T\cc\to TP$, $\Phi:T\cc\to TM$, which
are also defined on vectors based at the vertices of the tiling
$\cc=\cup_{g\in G}P_g$. Let $v$ be a vertex, and let $\al$ be the
angle of the corner $f(v)\in K(P)$. Then $\Phi:T_v\cc\to
T_{\vp(v)}M$ is $m$-to-$1$ if $\al=m\pi/n$ and $\infty$-to-$1$ if
$\al$ is $\pi$-irrational. The geodesics $\ga(t)$ in $\cc$ cannot
be further extended (generally) once they reach a vertex. All
other geodesics in $\cc$ are defined for $-\infty<t<\infty$.

Using the inclusion $P\subset M$, we identify $TP$ with the subset
of $TM$ consisting of $M$-tangent vectors with base-points in $P$,
and directed inward. Any $v\in TP$ defines the billiard orbit in
$P$, $\be(t)=\exp_P(tv),0\le t,$  and the geodesic in $M$,
$\ga(t)=\exp_M(tv),0\le t$. They are related by the canonical {\em
unfolding of billiard orbits}. This is an inductive procedure
which replaces the consecutive reflections about the sides of $P$
by consecutive reflections of the ``latest billiard table" $g(P)$
about the appropriate side, yielding the next billiard table
$h(P)$, and continuing the geodesic straight across the common
side of $g(P)$ and $h(P)$. See \cite{Gut86} in the planar case and
\cite{GT04}, section 3.1, in the general case. Let $x\in P$ and
let $v\in T_xP$. We denote by $\be_v$ (resp. $\ga_v$) the billiard
orbit in $P$ (resp. the geodesic in $M$) that emanates from $x$ in
the direction $v$. The {\em unfolding operator},
$U_x:\be_v\mapsto\ga_v$,  preserves the parametrisations.
\begin{prop}  \label{straigh_prop}
Let $x\in P,v\in T_xP$. Identify $P$ and $P_e\subset\cc$ and let
$x\in P_e,v\in T_x\cc$ be the corresponding data. Then for
$t\in\R_+$ we have
\begin{equation}   \label{straigh_eq}
U_x(\exp_P(v,t))=\vp(\exp_{\cc}(v,t)).
\end{equation}
\begin{proof}
We will freely use the preceding discussion.  As $t\in\R_+$ goes
to infinity, $\exp_P(v,t)$ runs with the unit speed along a
billiard orbit in $P$. The curve $\exp_{\cc}(v,t)$ is the geodesic
in $\cc$ defined by the data $(x,v)$, and $\vp(\exp_{\cc}(v,t))$
is the geodesic in $M$ emanating from $x$ in the direction $v$.
The billiard orbit in $P$ and the geodesic in $M$ are related by
the unfolding operator.
\end{proof}
\end{prop}

For $x\in P$ let $E_xP=T_xP\times\R_+$ be the full tangent space
(or the full tangent cone) at $x$.  If $S\subset T_xP$ is a
segment, let $ES_xP=S\times\R_+$ be the corresponding subcone. We
use the analogous notation for $x\in\cc$ or $x\in M$. In polar
coordinates $(t,\tht)$ in $\RR$ the lebesgue measure on $E_xP$ is
given by the density $tdtd\tht$.
%We view $\exp^x_P$ as a mapping from $E_xP$ into $P$.
%
\begin{corol}  \label{pull_back_cor}
Let $x\in P\subset M$ be arbitrary, and let $\exp^x_P:E_xP\to P$
be the exponential mapping. The pull-back by $\exp^x_P$ of the
lebesgue measure on $P$ to $E_xP$ is the smooth measure with the
density $d\nu(t,\tht)$.

\noindent 1. When $M=\RR$, we have $d\nu=tdtd\tht$.

\noindent 2. When $M=\hyp$, we have $d\nu=\sinh tdtd\tht$.

\noindent 3. When $M=\unisph$, we have $d\nu=|\sin t|dtd\tht$.

\begin{proof}
By Proposition~\ref{straigh_prop}, the measure in question
coincides with the pullback to the tangent space $E_xM$ of the
riemannian measure on $M$ by the exponential map $E_xM\to M$. The
latter is well known.
\end{proof}
\end{corol}

\medskip

We point out that the preceding material has a billiard map
version. We will briefly discuss it now. Let
$\be(t)=(z(t),\tht(t)),\,t\in\R,$ be an orbit of the billiard
flow. We obtain the corresponding billiard map orbit
$\be_d(k),\,k\in\Z,$ by restricting $\be(t)$ to the consecutive
times $t_k$ such that $z(t_k)\in\bo P$. The correspondence
$\be(\cdot)\mapsto\be_d(\cdot)$ is invertible. This allows us to
formulate the billiard map versions of the universal covering
space, the lifting of billiard map orbits to the universal
covering space, and the relationship between the liftings and the
unfoldings, \`a l\`a Proposition~\ref{straigh_prop}.  Since we are
not directly using this material in the body of the paper, we
spare the details.

%\medskip

%\newpage

%\medskip

%\medskip

%

\end{document}